\documentclass[12pt]{amsart}

\usepackage{amsmath,xspace,amssymb,mathrsfs}
\usepackage{color}

\input xy
\xyoption{all}
\xyoption{2cell}
\UseAllTwocells
\CompileMatrices

\newcommand{\Spec}{\operatorname{Spec}}
\renewcommand{\phi}{\varphi}

\newcommand{\Ker}{\operatorname{Ker}}

\newcommand{\Ima}{\operatorname{Im}}

\newcommand{\MA}{\operatorname{Max}}

\newcommand{\Min}{\operatorname{Min}}
\newcommand{\Ann}{\operatorname{Ann}}

\newcommand{\Supp}{\operatorname{Supp}}

\newtheorem{proposition}{Proposition}[section]
\newtheorem{lemma}[proposition]{Lemma}

\newtheorem{corollary}[proposition]{Corollary}
\newtheorem{theorem}[proposition]{Theorem}

\theoremstyle{definition}

\newtheorem{remark}[proposition]{Remark}

\begin{document}

\title{Notes on finitely generated flat modules}

\author[A. Tarizadeh]{Abolfazl Tarizadeh}
\address{Department of Mathematics, Faculty of Basic Sciences, University of Maragheh \\
P. O. Box 55136-553, Maragheh, Iran.
 }
\email{ebulfez1978@gmail.com}

\date{}
\footnotetext{ 2010 Mathematics Subject Classification: 13C10, 19A13, 13C11, 13E99.
\\ Key words and phrases: flat module; flat topology; patch topology; projectivity; S-ring.}

\begin{abstract} In this article, the projectivity of finitely generated flat modules of a commutative ring are studied from a topological point of view. Then various interesting results are obtained. For instance, it is shown that if a ring has either a finitely many minimal primes or a finitely many maximal ideals then every finitely generated flat module over it is projective. It is also shown that if a particular subset of the prime spectrum of a ring satisfies some certain ascending or descending chain conditions then every finitely generated flat module over this ring is projective. These results generalize some major results in the literature on the projectivity of finitely generated flat modules.
\end{abstract}

\maketitle

\section{Introduction}

Studying the projectivity of finitely generated flat modules has been the main topic of many articles over the years and it is still of current interest, see e.g. \cite{Cox-Pendleton}, \cite{Endo}, \cite{Jondrup}, \cite[\S4E]{Lam}, \cite{Puninski-Rothmaler}, \cite{Vasconcelos} and \cite{Wiegand}. The main motivation behind in the investigating the projectivity of f.g. flat modules stems from the fact that ``every f.g. flat module over a local ring is free''. We use f.g. in place of ``finitely generated''. Note that in general there are f.g. flat modules which are not necessarily projective, see \cite[Example 3.11]{Abolfazl} see also \cite[Tag 00NY]{Johan}. \\

It this article we have successfully applied the spectral (Zariski and flat) and patch topologies of the prime spectrum $\Spec R$ in order to investigate the projectivity of f.g. flat $R-$modules. The obtained results from this method generalize some major results in the literature on the projectivity of f.g. flat modules. In fact, Theorem \ref{coro 112} vastly generalizes \cite[Theorem 4.38]{Lam}, \cite[Corollary 1.5]{Jondrup},
\cite[Fact 7.5]{Puninski-Rothmaler} and \cite[Corollary 3.57]{Rotman} in the commutative case. Also Theorem \ref{th 55} generalizes \cite[Proposition 7.6]{Puninski-Rothmaler}. \\

Theorems \ref{th 2}, \ref{Theorem I}, \ref{coro 112} and \ref{th 55} and Corollaries \ref{Corollary I}, \ref{Corollary II} and \ref{Corollary III} are the main results of this article. In this article, all of the rings are commutative. \\

\section{Preliminaries}

Let $R$ be a commutative ring. Then there is a (unique) topology over $\Spec(R)$ such that the collection of subsets $V(f)=\{\mathfrak{p}\in\Spec R: f\in\mathfrak{p}\}$ with $f\in R$ formes a sub-basis for the opens of  this topology. It is called flat (or, inverse) topology.
Therefore, the collection of subsets $V(I)$ where $I$ runs through the set of finitely generated ideals of $R$ forms a basis for the flat opens. It is proved that the flat closed subsets of $\Spec(R)$ are precisely of the form $\Ima\phi^{\ast}$ where $\phi:R\rightarrow A$ is a ``flat'' ring map. Moreover there is a (unique) topology over $\Spec(R)$ such that the collection of subsets $D(f)\cap V(g)$ with $f,g\in R$ formes a sub-basis for the opens of this topology. It is called the patch (or, constructible) topology. Therefore the collection of subsets $D(f)\cap V(I)$ with $f\in R$ and $I$ runs through the set of finitely generated ideals of $R$ is a basis for the patch opens. It is also proved that the patch closed subsets of $\Spec(R)$ are precisely of the form $\Ima\phi^{\ast}$ where $\phi:R\rightarrow A$ is a ring map. Clearly the patch topology is finer than the flat and Zariski topologies. The flat topology behaves as the dual of the Zariski topology. For instance, if $\mathfrak{p}$ is a prime ideal of $R$ then its closure with respect to the flat topology comes from the canonical ring map $R\rightarrow R_{\mathfrak{p}}$. In fact, $\Lambda(\mathfrak{p})=\{\mathfrak{q}\in\Spec R: \mathfrak{q}\subseteq\mathfrak{p}\}$. Here $\Lambda(\mathfrak{p})$ denotes the closure of $\{\mathfrak{p}\}$ in $\Spec R$ with respect to the flat topology. Recall that a subset $E$ of $\Spec(R)$ is said to be stable under the generalization (resp. specialization) if  for any two prime ideals $\mathfrak{p}$ and $\mathfrak{q}$ of $R$ with $\mathfrak{p}\subset\mathfrak{q}$ (resp. $\mathfrak{q}\subset\mathfrak{p}$) if $\mathfrak{q}\in E$ then $\mathfrak{p}\in E$. One can show that a subset of $\Spec R$ is flat closed if and only if it is patch closed and stable under the generalization. Dually, a subset of $\Spec R$ is Zariski closed if and only if it is patch closed and stable under the specialization. A subset of $\Spec(R)$ is said to be a double-closed if it is closed with respect to the both flat and Zariski topologies. For more details see \cite{Hochster} or \cite{Abolfazl 2}. We shall freely use the above facts in this article. \\

\section{Main results}

In a topological space, the arbitrary unions of closed subsets are not necessarily closed. But in the prime spectrum, under some circumstances, interesting things are happening. To realize this, we first need to introduce two notations. Let $R$ be a ring. If $E$ is a subset of $\Spec(R)$ then we define $\mathcal{F}(E)=\bigcup\limits_{\mathfrak{p}\in E}\Lambda(\mathfrak{p})$ and $\mathcal{Z}(E)=\bigcup\limits_{\mathfrak{p}\in E}V(\mathfrak{p})$. We have then the following interesting result. \\

\begin{theorem}\label{th 2}
$\mathbf{(i)}$ If $E$ is Zariski closed then $\mathcal{F}(E)$ is flat closed.\\
$\mathbf{(ii)}$ If $E$ is patch closed then $\mathcal{Z}(E)$ is Zariski closed.\\
$\mathbf{(iii)}$ If $E$ is stable under the generalization and $E=V(I)$ for some ideal $I$ of $R$ then $R/J$ is $R-$flat and $E=V(J)$ where $J$ is the kernel of the canonical map $R\rightarrow S^{-1}R$ with $S=1+I$. \\
\end{theorem}

{\bf Proof.} $\textbf{(i)}$: Suppose $E=V(I)$ for some ideal $I$ of $R$. We claim that $\mathcal{F}(E)=\Ima\pi^{\ast}$ where $\pi:R\rightarrow S^{-1}R$ is the canonical map with $S=1+I$. The inclusion $\mathcal{F}(E)\subseteq\Ima\pi^{\ast}$ is obvious. To prove the reverse inclusion, let $\mathfrak{q}$ be a prime ideal of $R$ such that $\mathfrak{q}\cap S=\emptyset$. There exists a prime ideal $\mathfrak{p}$ of $R$ such that $\mathfrak{q}\subseteq\mathfrak{p}$ and $S^{-1}\mathfrak{p}$ is a maximal ideal of $S^{-1}R$. We have $I\subseteq\mathfrak{p}$. If not, then choose some element $f\in I\setminus\mathfrak{p}$. Clearly $(\mathfrak{p}+Rf)\cap S\neq\emptyset$.
Thus there are elements $r\in R$ and $g\in I$ such that $1+rf+g\in\mathfrak{p}$. But this is a contradiction since $\mathfrak{p}\cap S=\emptyset$. Therefore $\mathfrak{q}\in\mathcal{F}(E)$.\\
$\textbf{(ii)}$: We have $E=\Ima\phi^{\ast}$ for some ring morphism $\phi:R\rightarrow A$. It follows that $\mathcal{Z}(E)=V(I)$ where $I=\Ker\phi$. Because the inclusion $\mathcal{Z}(E)\subseteq V(I)$ is obvious. To prove the reverse inclusion, pick
$\mathfrak{p}\in V(I)$. Let $\mathfrak{q}$ be a minimal prime of $I$ such that $\mathfrak{q}\subseteq\mathfrak{p}$. Thus there exists a prime ideal of $A$ which lying over $\mathfrak{q}/I$. It follows that
$\mathfrak{q}\in E$. \\
$\textbf{(iii)}$: Clearly $E\subseteq V(J)$. By the proofs of $(i)$ and $(ii)$, we have $\mathcal{Z}\big(\mathcal{F}(E)\big)=V(J)$. It follows that $V(J)\subseteq E$ because $E$ is stable under the generalization. Thus $E=V(J)$. Using this, then by \cite[Theorem 3.8]{Abolfazl}, $R/J$ is $R-$flat.  $\Box$ \\

As a first application of Theorem \ref{th 2} we get the following result. \\

\begin{theorem}\label{Theorem I} Let $R$ be a ring. Then the assignment $I\rightsquigarrow V(I)$ is a bijective map from the set of ideals $I$ of $R$ such that $R/I$ is $R-$flat onto the set of Zariski closed subsets of $\Spec(R)$ which are stable under the generalization. \\
\end{theorem}

{\bf Proof.} First we show that this map is well-defined. That is, if $R/I$ is $R-$flat then we have to show that $V(I)$ is stable under the generalization. Let $\mathfrak{p}$ and $\mathfrak{q}$ be two prime ideals of $R$ such that $I\subseteq\mathfrak{p}$ and $\mathfrak{q}\subseteq\mathfrak{p}$. Suppose there is some $f\in I$ such that $f\notin\mathfrak{q}$. It follows that $\Ann(f)\subseteq\mathfrak{q}$. So $\Ann(f)+I\subseteq\mathfrak{p}$. But this is a contradiction since $\Ann(f)+I=R$. Thus $\mathfrak{q}\in V(I)$. Then we show that this map is injective. Let $I$ and $J$ be two ideals of $R$ such that $R/I$ and $R/J$ are $R-$flat and $V(I)=V(J)$. Take $f\in I$. If $f\notin J$ then by \cite[Corollary 3.9]{Abolfazl}, $\Ann(f)+J\neq R$. Thus there exists a prime ideal $\mathfrak{p}$ of $R$ such that $\Ann(f)+J\subseteq\mathfrak{p}$. It follows that $\Ann(f)+I\subseteq\mathfrak{p}$. But this is a contradiction since $\Ann(f)+I=R$. Therefore $I=J$. The surjectivity of this map implies from Theorem \ref{th 2} (iii). $\Box$ \\

\begin{remark} Note that a subset of $\Spec(R)$ is Zariski closed and stable under the generalization if and only if it is flat closed and stable under the specialization, see \cite[Theorem 3.11]{Abolfazl 2}. \\
\end{remark}

\begin{corollary}\label{Corollary I} Let $I$ be an ideal of a ring $R$ such that $R/\sqrt{I}$ is $R-$flat. Then $I=\sqrt{I}$. \\
\end{corollary}

{\bf Proof.} If $f\in I$ then by \cite[Theorem 3.8]{Abolfazl},  $\Ann(f)+\sqrt{I}=R$. It follows that $\sqrt{\Ann(f)}+\sqrt{I}=R$ and so $\Ann(f)+I=R$. Thus again by \cite[Theorem 3.8]{Abolfazl}, $R/I$ is $R-$flat. Then, by Theorem \ref{Theorem I}, $I=\sqrt{I}$. $\Box$ \\

\begin{corollary}\label{Corollary II} Let $I$ be an ideal of a reduced ring $R$ such that $R/I$ is $R-$flat. Then $I=\sqrt{I}$. \\
\end{corollary}

{\bf Proof.} By \cite[Corollary 3.9]{Abolfazl}, $\Supp(I)=\Spec(R)\setminus V(I)\subseteq\Supp(\sqrt{I})$. Conversely, if $\mathfrak{p}\in\Supp(\sqrt{I})$ then there exists some $f\in\sqrt{I}$ such that $f/1\neq0$. If $I_{\mathfrak{p}}=0$ then there exist $s\in R\setminus\mathfrak{p}$ and a natural number $n\geq1$ such that $sf^{n}=0$. It follows that $sf=0$ since $R$ is reduced. But this is a contradiction. Therefore $\mathfrak{p}\in\Supp(I)$. Hence, $R/\sqrt{I}$ is $R-$flat. Thus, by Theorem \ref{Theorem I}, $I=\sqrt{I}$. $\Box$ \\

\begin{corollary}\label{Corollary III} Let $I$ and $J$ be two ideals of a reduced ring $R$ such that $R/I$ is $R-$flat and $V(I)=V(J)$. Then $I=J$. \\
\end{corollary}

{\bf Proof.} It is an immediate consequence of Corollaries \ref{Corollary I} and \ref{Corollary II}. $\Box$ \\

\begin{lemma}\label{lemma 9} Let $I$ be an ideal of a ring $R$. If $R/I$ is $R-$flat then for each finite subset $\{f_{1},...,f_{n}\}$ of $I$ there exists some $g\in I$ such that $f_{i}=f_{i}g$ for all $i$.\\
\end{lemma}

{\bf Proof.} By \cite[Theorem 3.8]{Abolfazl}, $R/I$ is $R-$flat if and only if  $\Ann_{R}(f)+I=R$ for all $f\in I$. Thus for each pair $(f,f')$ of elements of $I$ then there exist $h,h'\in I$ such that $f=fh$ and $f'=f'h'$. Clearly $g:=h+h'-hh'\in I$, $f=fg$ and $f'=f'g$. $\Box$ \\

A ring $R$ is called an S-ring (``S'' refers to Sakhajev) if every f.g. flat $R-$module is $R-$projective. \\

We have improved the following result by adding $\mathbf{(iv)-(vii)}$ as new equivalents. The equivalency of the classical criteria are also proved by new methods. \\

\begin{theorem}\label{th 1} For a ring $R$ the following conditions are equivalent.\\
$\mathbf{(i)}$ The ring $R$ is an S-ring.\\
$\mathbf{(ii)}$ Every cyclic flat $R-$module is $R-$projective.\\
$\mathbf{(iii)}$ $R/I$ is $R-$projective whenever it is $R-$flat where $I$ is an ideal of $R$.\\
$\mathbf{(iv)}$ Every Zariski closed subset of $\Spec(R)$ which is stable under the generalization is Zariski open.\\
$\mathbf{(v)}$ Every patch closed subset of $\Spec(R)$ which is stable under the generalization and specialization is patch open.\\
$\mathbf{(vi)}$ Every flat closed subset of $\Spec(R)$ which is stable under the specialization is flat open.\\
$\mathbf{(vii)}$ Each double-closed subset of $\Spec(R)$ is of the form $V(e)$ where $e\in R$ is an idempotent.\\
$\mathbf{(viii)}$ For every sequence $(f_{n})_{n\geq1}$ of elements of $R$ if $f_{n}=f_{n}f_{n+1}$ for all $n$ then there exists some $k$ such that $f_{k}$ is an idempotent and $f_{n}=f_{k}$ for all $n\geq k$.\\
$\mathbf{(ix)}$ For every sequence $(g_{n})_{n\geq1}$ of elements of $R$ if $g_{n+1}=g_{n}g_{n+1}$ for all $n$ then there exists some $k$ such that $g_{k}$ is an idempotent and $g_{n}=g_{k}$ for all $n\geq k$.\\
\end{theorem}

{\bf Proof.} The implications $\textbf{(i)}\Rightarrow\textbf{(ii)}\Rightarrow\textbf{(iii)}$ are obvious. \\
$\textbf{(iii)}\Rightarrow\textbf{(iv)}:$ Suppose $E\subseteq\Spec(R)$ is stable under the generalization and $E=V(I)$ for some ideal $I$ of $R$. By Theorem \ref{th 2}, there is an ideal $J$ such that
$R/J$ is $R-$flat and $E=V(J)$. Thus, by the hypothesis, $R/J$ is $R-$projective. It follows that $E$ is Zariski open because the support of a projective module is Zariski open. \\
$\textbf{(iv)}\Leftrightarrow\textbf{(v)}\Leftrightarrow\textbf{(vi)}:$ See the Introduction. \\
$\textbf{(vi)}\Rightarrow\textbf{(vii)}:$ Apply \cite[Tag 00EE]{Johan}. \\
$\textbf{(vii)}\Rightarrow\textbf{(iv)}:$ There is nothing to prove. \\
$\textbf{(iv)}\Rightarrow\textbf{(i)}:$ Let $M$ be a f.g. flat $R-$module. To prove the assertion, by \cite[Theorem 3.2]{Abolfazl},
it suffices to show that for each natural number $n$, $\psi^{-1}(\{n\})$ is Zariski open where $\psi$ is the rank map of $M$, see \cite[Remark 2.7]{Abolfazl}. We have $\psi^{-1}(\{n\})=\Supp N\cap\big(\Spec(R)\setminus\Supp N'\big)$ where $N=\Lambda^{n}(M)$ and $N'=\Lambda^{n+1}(M)$. But $\Supp N$ and $\Supp N'$ are Zariski closed since $N$ and $N'$ are f.g. $R-$modules. But $N$ is a flat $R-$module. By applying \cite[Theorem 2.3]{Abolfazl} then we observe that the support of a f.g. flat module is stable under the generalization. Thus by the hypothesis, $\Supp N$ is Zariski open. Therefore $\psi^{-1}(\{n\})$ is Zariski open.\\
$\textbf{(iii)}\Rightarrow\textbf{(viii)}:$ Let $I=(f_{n}: n\geq1)$.
Clearly $\Ann_{R}(f)+I=R$ for all $f\in I$. It follows that
$R/I$ is $R-$flat and so, by the hypothesis, it is $R-$projective. Therefore by \cite[Lemma 3.1]{Abolfazl}, there exists $g\in I$ such that $I=Rg$. It follows that there is some $d\geq1$ such that $Rg=Rf_{d}$ since $I=\bigcup\limits_{n\geq1}Rf_{n}$. Let $k=d+1$. There exists some $r\in R$ such that $f_{k}=rf_{d}=rf_{d}f_{k}=f^{2}_{k}$. We also have $f_{k+1}=r'f_{k}$
for some $r'\in R$. It follows that $f_{k+1}=f_{k+1}f_{k}=f_{k}$ and by the induction we obtain that $f_{n}=f_{k}$ for all $n\geq k$.\\
$\textbf{(viii)}\Rightarrow\textbf{(iii)}:$ Let $I$ be an ideal of $R$ such that $R/I$ is $R-$flat. We shall prove that $I$ is generated by an idempotent element. To do this we act as follows. Let $\mathcal{I}$ be the set of ideals of the form $Re$ where $e\in I$ is an idempotent element. Let $\{Re_{n}: n\geq1\}$ be an ascending chain of elements of $\mathcal{I}$. For each $n$ there is some $r_{n}\in R$ such that $e_{n}=r_{n}e_{n+1}$. It follows that $e_{n}=e_{n}e_{n+1}$. Thus, by the hypothesis, the chain $Re_{1}\subseteq Re_{2}\subseteq ...$ is stationary. Therefore, by the axiom of choice, $\mathcal{I}$ has at least a maximal element $Re$. We also claim that if $J=(f_{n}: n\geq1)$ is a countably generated ideal of $R$ with $J\subseteq I$ then there exists an idempotent $e'\in I$ such that $J\subseteq Re'$. Because, by Lemma \ref{lemma 9}, there is an $g_{1}\in I$ such that $f_{1}=f_{1}g_{1}$. Then for the pair $(g_{1},f_{2})$, again by Lemma \ref{lemma 9}, we may find an $g_{2}\in I$ such that $g_{1}=g_{1}g_{2}$ and $f_{2}=f_{2}g_{2}$. Therefore, in this way, we obtain a sequence $(g_{n})$ of elements of $I$ such that $J\subseteq L=(g_{n} : n\geq1)$ and $g_{n}=g_{n}g_{n+1}$ for all $n\geq1$. But, by the hypothesis, there exists some $k\geq 1$ such that $g_{k}$ is an idempotent and $g_{n}=g_{k}$ for all $n\geq k$. It follows that $L=Rg_{k}$. This establishes the claim. Now pick $f\in I$. Then, by what we have proved above, there is an idempotent $e'\in I$ such that $Re\subseteq(e,f)\subseteq Re'$. By the maximality of $Re$, we obtain that $e=e'$. Thus $I=Re$ and so $R/I$ as $R-$module is isomorphic to $R(1-e)$. Therefore $R/I$ is $R-$projective.\\
$\textbf{(viii)}\Leftrightarrow\textbf{(ix)}:$ Let $(f_{n})$ be a sequence of elements of $R$. Put $g_{n}:=1-f_{n}$ for all $n$. Then $f_{n}=f_{n}f_{n+1}$ if and only if $g_{n+1}=g_{n}g_{n+1}$. $\Box$ \\

As a consequence of Theorem \ref{th 1}, we obtain the following result which in turn vastly generalizes some previous results in the literature specially including \cite[Theorem 4.38]{Lam}, \cite[Corollary 1.5]{Jondrup},
\cite[Fact 7.5]{Puninski-Rothmaler} and \cite[Corollary 3.57]{Rotman} in the commutative case.\\

\begin{theorem}\label{coro 112} Let $R$ be a ring which has either a finitely many minimal primes or a finitely many maximal ideals. Then $R$ is an S-ring.\\
\end{theorem}

{\bf Proof.} Let $F$ be a patch closed subset of $\Spec(R)$ which is stable under the generalization and specialization. By Theorem \ref{th 1}, it suffices to show that it is a patch open. First assume that $\Min(R)=\{\mathfrak{p}_{1},...,\mathfrak{p}_{n}\}$.
There exists some $s$ with $1\leq s\leq n$ such that $\mathfrak{p}_{s},\mathfrak{p}_{s+1},...,\mathfrak{p}_{n}\notin F$ but $\mathfrak{p}_{i}\in F$ for all $i<s$.
It follows that $\Spec(R)\setminus F=\bigcup\limits_{i=s}^{n}V(\mathfrak{p}_{i})$. Therefore $F$ is Zarsiki open in this case and so it is patch open. Now assume that $\MA(R)=\{\mathfrak{m}_{1},...,\mathfrak{m}_{d}\}$. Similarly, there exists some $k$ with $1\leq k\leq d$ such that $\mathfrak{m}_{k},\mathfrak{m}_{k+1},...,\mathfrak{m}_{d}\notin F$ but $\mathfrak{m}_{i}\in F$ for all $i<k$.
We have $\Spec(R)\setminus F=\bigcup\limits_{i=k}^{d}\Lambda(\mathfrak{m}_{i})$. Therefore $F$ is a flat open in this case and so it is patch open. $\Box$ \\

\begin{remark} In relation with Theorem \ref{coro 112}, note that though very projective module over a local ring is free (\cite[Tag 0593]{Johan}) but in general this is not necessarily true even for a semi-local ring (a ring with a finitely many maximal ideals) which is not local. As a specific example, let $n>1$ be a natural number which has at least two distinct prime factors and let $n=p^{s_{1}}_{1}...p^{s_{k}}_{k}$ be its prime factorization where the $p_{i}$ are distinct prime numbers and $s_{i}\geq1$ for all $i$. Each $A_{i}$ can be considered as $R-$module through the canonical ring map $R\rightarrow A_{i}$ where $R=\mathbb{Z}/n\mathbb{Z}$ and $A_{i}=\mathbb{Z}/p^{s_{i}}_{i}\mathbb{Z}$. By the Chinese remainder theorem, $R$ as module over itself is isomorphic to the direct sum $\bigoplus\limits_{i=1}^{k}A_{i}$. Thus each $A_{i}$ is $R-$projective. But none of them is $R-$free since every non-zero free $R-$module has at least $n$ elements while $p^{s_{i}}_{i}<n$ for all $i$. Note that $R$ is a semi-local ring with the maximal ideals $p_{i}\mathbb{Z}/n\mathbb{Z}$.\\
\end{remark}

\begin{theorem}\label{th 55} Let $X$ be a subset of $\Spec(R)$ with the property that for each maximal ideal $\mathfrak{m}$ of $R$ there exists some $\mathfrak{p}\in X$ such that $\mathfrak{p}\subseteq\mathfrak{m}$. If the collection of subsets $X\cap V(f)$ with $f\in R$ satisfies either the ascending chain condition or the descending chain condition then $R$ is an S-ring. \\
\end{theorem}

{\bf Proof.} By \cite[Theorem 3.5]{Abolfazl}, it suffices to show that $R/J$ is an S-ring where $J=\bigcap\limits_{\mathfrak{p}\in X}\mathfrak{p}$.
Let $(x_{n})$ be a sequence of elements of $R/J$ such that $x_{n}=x_{n}x_{n+1}$ for all $n$. Suppose $x_{n}=a_{n}+J$ for all $n$.
Let $E_{n}=X\cap V(a_{n})$ and let $F_{n}=X\cap V(1-a_{n})$. Clearly $E_{n}\supseteq E_{n+1}$, $F_{n}\subseteq F_{n+1}$ and $X=E_{n}\cup F_{n+1}$. First assume the descending chain condition. Then there exists some $d\geq1$ such that $E_{n}=E_{d}$ for all $n\geq d$. Therefore $X=E_{n}\cup F_{n}$ for all $n>d$. Thus $a_{n}(1-a_{n})\in\mathfrak{p}$ for all  $\mathfrak{p}\in X$ and all $n>d$. It follows that $x_{n}=x^{2}_{n}$ for all $n>d$. The chain of ideals $(x_{d+1})\subseteq(x_{d+2})\subseteq...$ eventually stabilizes. If not, then the ascending chain $V(1-x_{d+1})\subseteq V(1-x_{d+2})\subseteq...$ does not stabilize. Therefore we may find some $k>d$ such that $V(1-x_{k})$ is a proper subset of $V(1-x_{k+1})$. Thus there exists a prime ideal $\mathfrak{q}$ of $R$ such that $J\subseteq\mathfrak{q}$ and
$1-a_{k+1}, a_{k}\in\mathfrak{q}$. There is a maximal ideal $\mathfrak{m}$ of $R$ such that $\mathfrak{q}\subseteq\mathfrak{m}$. By the hypotheses, there is a $\mathfrak{p}\in X$ such that $\mathfrak{p}\subseteq\mathfrak{m}$. Clearly  $1-a_{k+1}, a_{k}\in\mathfrak{p}$. This means that $E_{k+1}$ is a proper subset of $E_{k}$. But this is a contradiction. Thus there is some $s>d$ such that $(x_{n})=(x_{s})$ for all $n\geq s$. Therefore $x_{n}=x_{s}$ for all $n\geq s$ and so by Theorem \ref{th 1}, $R/J$ is an S-ring in the case of the descending chain condition. Apply a similar argument as above for the chain $F_{1}\subseteq F_{2}\subseteq...$ in the case of the ascending chain condition.  $\Box$ \\

\begin{corollary}\label{th 5} If the collection of subsets $\Min(R)\cap V(f)$ with $f\in R$ satisfies either the ascending chain condition or the descending chain condition then $R$ is an S-ring. \\
\end{corollary}

{\bf Proof.} It implies from Theorem \ref{th 55} by taking $X=\Min(R)$. $\Box$ \\

\begin{corollary}\cite[Proposition 7.6]{Puninski-Rothmaler}\label{th 515} If the collection of subsets $\MA(R)\cap V(f)$ with $f\in R$ satisfies either the ascending chain condition or the descending chain condition then $R$ is an S-ring. \\
\end{corollary}

{\bf Proof.} In Theorem \ref{th 55}, put $X=\MA(R)$. $\Box$ \\

\begin{proposition} The direct product of a family of rings $(R_{i})_{i\in I}$ is an S-ring if and only if $I$ is a finite set and each $R_{i}$ is an S-rings.\\
\end{proposition}

{\bf Proof.} Let $R=\prod\limits_{i\in I}R_{i}$ be an S-ring. We may assume that all of the rings $R_{i}$ are non-zero. Suppose $I$ is an infinite set. Consider a well-ordered relation $<$ on $I$. Let $i_{1}$ be the lest element of $I$ and for each natural number $n\geq1$, by induction, let $i_{n+1}$ be the least element of $I\setminus\{i_{1},...,i_{n}\}$. Now we define $x_{n}=(r_{n,i})_{i\in I}$ as an element of $R$ by $r_{n,i}=1$ for all $i\in\{i_{1},...,i_{n}\}$ and $r_{n,i}=0$ for all $i\in I\setminus\{i_{1},...,i_{n}\}$.
Clearly the sequence $(x_{n})$ satisfies the condition $x_{n}=x_{n}x_{n+1}$. Thus, by Theorem \ref{th 1}, there is some $k$ such that $x_{n}=x_{k}$ for all $n\geq k$. But this is a contradiction. Thus $I$ should be a finite set.
The remaining assertions, by applying Theorem \ref{th 1}(viii), are straightforward. $\Box$ \\

\end{document}